\documentclass[a4paper,12pt]{article}

\usepackage{amsmath, wrapfig}
\usepackage{dsfont, a4wide, amsthm, amssymb, amsfonts, graphicx}
\usepackage{fancyhdr, xspace, psfrag, setspace, supertabular, color, enumerate}
\usepackage{hyperref}
\usepackage{bbm}

\newtheorem{theorem}{Theorem}[section]
\newtheorem{proposition}[theorem]{Proposition}
\newtheorem{lemma}[theorem]{Lemma}

\newtheorem{definition}[theorem]{Definition}

\newtheorem{example}[theorem]{Example}

\newtheorem*{example*}{Example}

\def\P{\mathbb{P}}
\def\R{\mathbb{R}}

\def\N{\mathbb{N}}

\title{Essential and disjoint edge numbers for graphs and hypergraphs}
\author{Thomas Karam\footnote{School of Mathematical Sciences, Shanghai Jiao Tong University. Email: \texttt{thomaskaram@sjtu.edu.cn}. This research is supported by the National Natural Science Foundation of China, grant number W2633011.}}

\begin{document}
\maketitle

\begin{abstract}

We study the extension to arbitrary forbidden sets of the property that every graph has a bipartition over which at least half of the graph edges are cross-edges, and obtain analogous results in the hypergraph setting. 

\end{abstract}

\tableofcontents

\section{Introduction}

\subsection{Broader context: essential versus disjoint complexity}

Consider first a class $\mathcal{O}$ of multidimensional objects. That class may, for instance, consist of subsets of $\R^d$ or more generally functions defined on $\R^d$. Second, consider a notion of complexity on the class of objects. That notion of complexity may for instance be size, Lebesgue measure, a covering number with respect to some “simple” objects such as subspaces, matrix rank, a notion of rank on tensors and so forth. We now suppose that the object is unknown or modifiable on some subset $U$ of $\R^d$ and known everywhere else. 

This paper aims to contribute to developing a systematic understanding of a min-max phenomenon in which for suitable $U$, the smallest complexity that an object $A \in \mathcal{O}$ can have after any modifications inside $U$ and the largest complexity that can be witnessed by restricting $A$ to a Cartesian box avoiding $U$ are two sides of the same coin, in a way that is uniform with respect to the choice of $A \in \mathcal{O}$. This min-max phenomenon has in the past received attention for several classes of objects, in particular graphs and hypergraphs \cite{Conlon Fox Kwan Sudakov} \cite{Erdos} \cite{Erdos Kleitman} \cite{Janzer Portier} \cite{Poljak Tuza} \cite{Scott} as we are about to discuss in more detail in the next subsection, covering numbers of lattice subsets \cite{Karam simple}, and ranks of tensors \cite{Karam highrank}. Yet in these situations, the set $U$ has usually been a very specific type of set, such as the diagonal, and we here initiate a study emphasising uniformity with respect to the choice of the set $U$ (within a suitable class).

In this paper we focus on what appears to be the most basic setting, where the class $\mathcal{O}$ consists of finite subsets of $\N^d$ and the notion of complexity is the counting measure. In this setting, a major consideration that simplifies our analysis is the uniqueness of decompositions, in the following sense. The size of a subset $A$ of $\N^d$ is the smallest nonnegative integer $k$ such that we may write \[A \subset A_1 \cup \dots \cup A_k\] for singletons $A_1,\dots,A_k$. Provided that $k$ is minimal, the sets $A_1,\dots,A_k$ are uniquely determined up to permutation. In contrast, many notions of complexity share with the size the fact that the complexity of an object $A$ is defined in terms of the smallest number of objects $A_1,\dots,A_k$ with complexity one necessary to cover or decompose $A$ into, yet the choice of these objects $A_1,\dots,A_k$ with complexity one is far from uniquely determined. To take the example of covering numbers of finite sets of Euclidean space with lines, there is often more than one covering of a finite set of points with the minimum possible number of lines: this is easily seen by taking three non-aligned points and covering them with two lines. Likewise, for the tensor rank, smallest-length decompositions are typically not unique, with uniqueness of decompositions being itself a major topic of study. Still with the tensor rank (also called the CP-rank), this question of which tensors have a unique decomposition (up to permutation) is called the identifiability question, and numerous relevant references are discussed for instance in the paper \cite{Ballico Bernardi Santarsiero}.

As will be discussed in the next subsection, the case of set sizes arises naturally as a generalisation of a well-known property on graph partitioning (and often one of the first applications of the probabilistic method that one encounters). However, in the light of the more general context above we intend the present paper to also act as a reference point for later arguments arising in the case of complexity notions that are harder to study. Can the argument be adapted straightforwardly and if not, then where do the main complications arise, and what are the corresponding workarounds?

\subsection{Intrinsic context: graph and hypergraph restrictions}

If $G$ is a graph or hypergraph, then we write $e(G)$ for the number of edges of $G$.

One of the very first applications of the probabilistic method commonly encountered is the fact that for every finite undirected loopless graph $G = (V,E)$, there exists a bipartition $\{X,Y\}$ of the vertex set $V$ (allowing empty parts) such that the number $e(G[X,Y])$ of edges of $G$ joining a vertex of $X$ and a vertex of $Y$ is at least half of the total number $e(G)$ of edges of $G$. To establish this fact proved by Erd\H{o}s \cite{Erdos} in 1965, it suffices to send the vertices of $G$ independently and uniformly at random to the sets $X$ and $Y$; for every edge of $G$, the probability that that edge joins a vertex of $X$ and a vertex of $Y$ is equal to $1/2$, and it then follows (by linearity of expectation) that there exists in particular some satisfactory bipartition $\{X,Y\}$ of $V$.

Erd\H{o}s and Kleitman then extended both the result and its proof to uniform hypergraphs in 1968 \cite{Erdos Kleitman}. In the below proposition, a hyperedge of a $d$-uniform hypergraph $G=(V,E)$ is a set of $d$ pairwise distinct vertices of $V$.

\begin{proposition}\label{standard hypergraph bound}

Let $d \ge 2$ be an integer, let $G=(V,E)$ be a finite $d$-uniform undirected hypergraph. There exists a partition $\{X_1,\dots,X_d\}$ of $V$ (allowing empty parts), such that the proportion of hyperedges $\{x_1,\dots,x_d\}$ containing exactly one vertex in each part $X_1,\dots,X_d$ is at least $(d!/d^d) e(G)$.

\end{proposition}

Indeed, the proof from the graph case may be adapted, by sending every vertex of $V$ independently and uniformly at random to one of the sets $X_1,\dots,X_d$. The proportion $d!/d^d$ is then the proportion of satisfactory assignments, and may be seen to be optimal in Proposition \ref{standard hypergraph bound} by considering (large) complete $d$-uniform hypergraphs.

A survey on Proposition \ref{standard hypergraph bound} and related results was written by Scott \cite{Scott} in 2005. More recently still, progress has been achieved on the topics of hypergraph cuts, in particular by Conlon, Fox, Kwan and Sudakov \cite{Conlon Fox Kwan Sudakov} and by Janzer and Portier \cite{Janzer Portier}. 

This result has a directed analogue, with a proportion of $1/d^d$, and a proof that is in essentially all other respects identical. As that setting is slightly more general, and also lends itself better to comparison with various settings – such as coverings of points in Euclidean space, and ranks of matrices and tensors, for which coordinates may typically not be gratuitously swapped without additional symmetry assumptions – we will be working directed graphs and hypergraphs.

In the setting of the ranks of tensors, a question analogous to that directed analogue of Proposition \ref{standard hypergraph bound} has been considered by the author in \cite{Karam highrank}.

Already in the case of (henceforth, directed) graphs, we may ask for an analogous statement if the graph is allowed to contain loops. In that setting, the existence of the desired bipartition $\{X,Y\}$ may fail qualitatively even if $G$ contains many edges: the cause is that loops of $G$ never appear as edges in any of the induced subgraphs $G[X,Y]$ with $X,Y$ disjoint, and hence have no effect on the number of edges of these graphs. To rescue the statement in the case of graphs avoiding loops, it is for that reason necessary to first delete all loops of $G$ before expressing a useful lower bound on the maximum possible $e(G[X,Y])$ in terms of $e(G)$. Once these deletions have occurred, the resulting graph is loopless, which allows us to conclude immediately as follows. We write $\Delta$ for the set $\{(x,x): x \in V\}$ of all possible loops.

\begin{proposition}\label{graph with loops}

Let $G = (V,E)$ be a directed graph, possibly with loops. Then there exist disjoint subsets $X,Y$ of $V$ such that the number $e(G[X \times Y])$ of edges of $G$ from a vertex in $X$ to a vertex in $Y$ is at least $|E \setminus \Delta|/4$.

\end{proposition}

Requiring $X,Y$ to be disjoint is equivalent to requiring the product $X \times Y$ to avoid (that is, have empty intersection with) the diagonal $\Delta$. We may wish to generalise this statement to instead forbid any arbitrary collection of pairs in $X \times Y$. That is, we may wish to obtain, for arbitrary sets $U$ (inside $V \times V$) an analogue of Proposition \ref{graph with loops} which provides the existence of $X,Y$ with $(X \times Y) \cap U = \emptyset$ such that the induced subgraph $G[X \times Y]$ contains many edges.

In the same way that loops cannot have any effect on any induced subgraphs $G[X \times Y]$ with $X,Y$ disjoint, there are edges that cannot have any effect on any induced subgraphs $G[X \times Y]$ with $X \times Y$ avoiding $U$: these are precisely the edges in $U$. To guarantee the existence of the desired sets $X,Y$, it would therefore be necessary (but not necessarily sufficient) to count only the edges of $E$ that are outside $U$.

At this point, we may nonetheless still hope that there exists an absolute constant $c>0$ such that if $G$ is a directed graph (possibly with loops), then there exist $X,Y$ with $(X \times Y) \cap U = \emptyset$ and satisfying \[e(G[X,Y]) \ge c |E \setminus U|.\]

The following example however illustrates that such an inequality however is false in general, and moreover that the disjoint edge number does not control the essential edge number uniformly over all $U$.

\begin{example}\label{Complement of diagonal}

Let $V$ be a finite set, and let \[U = (V \times V) \setminus \Delta = \{(x,y) \in V \times V: x \neq y\}.\] We take $G  = (V,E)$ with $E = \Delta$. Then we have $|E \setminus U| = |\Delta| = |V|$, but if $X,Y \subset V$ are such that the product $X \times Y$ has empty intersection with $U$, then $X \times Y$ is contained in $\Delta$, so either one of the sets $X,Y$ is empty, or $X,Y$ are equal singletons, leading to $e(G[X,Y]) \le 1$. 

\end{example}

Let us introduce definitions and notations for the primary quantities that we have been discussing so far. Let $G=(V,E)$ be a directed graph (possibly with loops) and let $U$ be a subset of $V \times V$. We define the \emph{essential edge number of $G$ with respect to $U$}, denoted by $e_U(G)$, to be the number $|E \setminus U|$ of edges of $G$ outside $U$, and define the \emph{disjoint edge number of $G$ with respect to $U$}, denoted by $d_U(G)$ to be the maximum number $e(G[X,Y])$ over all pairs $(X,Y)$ of subsets of $V$ such that the product $X \times Y$ has empty intersection with $U$. For instance, Proposition \ref{graph with loops} states that $d_{\Delta}(G) \ge e_\Delta(G)/4$ for all $G$.

For any $G=(V,E)$ and any $U$, the disjoint edge number is always at most the essential edge number. Indeed, if $X,Y \subset V$ are such that $(X \times Y) \cap U = \emptyset$, then the graph $G[X,Y]$ is an induced subgraph of the graph $(V, E \setminus U)$, of which $e_U(G)$ is the number of edges. Whether and how strongly an inequality holds in the converse direction is the much more delicate question.

Choosing a subset $U \subset V \times V$ uniformly at random with $|V|$ large establishes that $d_U(G)$ may be considerably lower than $e_U(G)$ in general. Indeed, such a set $U$ will typically have size close to $|V|^2/2$, so taking $E = V \times V$ leads to $e_U(G)$ close to $|V|^2/2$, but a probabilistic calculation shows that growing a Cartesian box avoiding $U$ to contain remotely that many points will typically be impossible.

Our first result will be to identify a complexity notion on subsets of the type $U \subset V \times V$ in terms of which one can characterise whether the disjoint edge number with respect to $U$ controls the essential edge number with respect to $U$. We will now define this notion.

\begin{definition}

Let $V$ be a finite set and let $U \subset V \times V$. We say that an \emph{isolated set with respect to $U$} is a set of points $(x_1,y_1), \dots, (x_r,y_r)$ in the complement of $U$ in $V \times V$, such that for every $1 \le i<j \le r$ we have $(x_i,y_j) \in U$ or $(x_j,y_i) \in U$. We say that the \emph{isolation number of $U$} is the maximum size of isolation sets with respect to $U$. (We note that the empty set is an isolation set.)

\end{definition}

We note that in the above definition, the values $x_1,\dots,x_r$ must be pairwise distinct, and likewise the values $y_1,\dots,y_r$ must be pairwise distinct (as it is straightforward to check that if they are not, then the property cannot be satisfied).

With sets $U$ with a high isolation number, the number $d_U$ does not control the number $e_U$ well. Indeed, taking the edge set of $G$ to be an isolated set with respect to $U$ that has maximal size $r \ge 1$, it follows immediately from the definition of $e_U(G)$ that $e_U(G) = r$ (as all edges of $E$ are outside $U$). Meanwhile, no pair $(X,Y)$ of subsets of $V$ such that $X \times Y$ avoids $U$ can have both its sets $X,Y$ with size $2$ or more, so $d_U(G) = 1$.

\subsection{Main results bounding essential edge number in terms of disjoint edge number}

Our first main result shall be the following qualitative converse, which informally states that uniformly over all sets $U$ with isolation number at most $r$ (and uniformly over all graphs $G = (V,E)$) the disjoint edge number $d_U$ of $G$ controls its essential edge number $e_U$. For positive integers $k,l$ we write $R(k,l)$ for the $(k,l)$-th Ramsey number, that is, for the value of the smallest integer $n$ such that every complete graph with $n$ vertices and with edges each coloured in red or blue contains a red complete graph with $k$ vertices or a blue complete graph with $l$ vertices.

\begin{theorem}\label{Main theorem in the two-set case}

Let $V$ be a finite set, let $U$ be a subset of $V \times V$, and let $r$ be the isolation number of $U$. Then every graph $G = (V,E)$ with vertex set $V$ satisfies \[e_U(G) \le R(d_U(G)+1,r+1) \le (d_U(G)+r+1)^{r+1}.\]

\end{theorem}

In the above theorem, the second inequality follows immediately from the well-known Erd\H{o}s-Szekeres bound. As we will show, the resulting exponential bound is sharp up to the linear constant in $r$ in the exponent.
\begin{proposition} \label{exponential lower bound}

Let $r$ be an integer. Then for every integer $n$ large enough (depending on $r$) there exists a graph $G=(V,E)$ (with $|V|=n$) and a subset $U$ of $V \times V$ with isolation number at most $r$ such that $e_U(A) = n$ but $d_U(A) \le n^{5/r}$.

\end{proposition}

The proof of Theorem \ref{Main theorem in the two-set case} will involve, as its primary tool, an auxiliary graph with vertex set $E \setminus U$ and an application of Ramsey’s theorem. In summary, monochromatic cliques of one colour will correspond to isolated sets with respect to $U$ and monochromatic cliques of the other colour will correspond to sets of edges arising as edge sets of the induced subgraphs over which the maximum number is taken in the definition of the $d_U(G)$.

We next describe how these results extend to the setting of $d$-uniform directed hypergraphs. 

Let $d \ge 2$ be an integer, let $H=(V,E)$ be a $d$-uniform directed hypergraph, and let $U$ be a subset of $V^d$. We say that the \emph{essential edge number of $G$ with respect to $U$}, denoted by $e_U(H)$, is the size of the set $E \setminus U$, and that the \emph{disjoint edge number of $G$ with respect to $U$}, denoted by $d_U(H)$, is the largest size of an intersection $E \cap (X_1 \times \dots \times X_d)$, over $d$-tuples of subsets $(X_1,\dots,X_d)$ such that the Cartesian box $X_1 \times \dots \times X_d$ has empty intersection with $U$.

To formulate our characterisation, we first define an appropriate analogue of isolated sets to higher dimensions.

Let $d \ge 2$ be an integer, let $V$ be a finite set, and let $U \subset V^d$. Whenever $A$ is a subset of $V^d$, we write $\mathcal{E}(A)$ for the \emph{enclosure of $A$}, that is, the smallest Cartesian box enclosing $A$. The set $\mathcal{E}(A)$ is the product $\prod_{j=1}^d \pi_j(A)$ where $\pi_j$ is the projection on the $j$th coordinate axis. We say that an \emph{isolated set with respect to $U$} is a set of points $x_1,\dots,x_r$ in the complement of $U$ inside $V^d$, such that $\mathcal{E}(\{x_i,x_j\})$ has non-empty intersection with $U$ whenever $i,j$ are distinct elements of $[r]$. As in the graph case, the \emph{isolation number} of $U$ is the largest size of an isolated set with respect to $U$.

In the hypergraph case, the formulation of our main theorem is as follows.

\begin{theorem} \label{Main theorem for hypergraphs}
Let $d \ge 2$ be an integer, let $V$ be a finite set, let $U$ be a subset of $V^d$, and let $r$ be the isolation number of $U$. Then every $d$-uniform directed hypergraph $H = (V,E)$ satisfies $e_U(H) \le (2d_U(H))^{d2^dr+1}$.
\end{theorem}

We note that for a fixed power of $d$, our upper bounds remain linear-exponential with respect to $r$. Compared to the graph case, the proof will be significantly more involved. The key difficulty jump is the following. In the $d=2$ case, an element in the enclosure of a set $E$ inherits its coordinates from at most two elements of $E$, and stating that the enclosure of a set $E$ has empty intersection with $U$ is equivalent to stating that the enclosures of its pairs of points all do. In contrast, this equivalence no longer persists with higher values of $d$, as an element of in the enclosure of a set $E$ may inherit its coordinates from up to $d$ elements of $E$. Yet despite that, the parameter $r$, which is built from considerations only involving pairs of points, still controls the maximum possible growth rate of the essential edge number in terms of the disjoint edge number. In the proof, we will resort to building a tree that describes the structure of the elements of $E$ by keeping track of the subsets of coordinates that change.

In Section \ref{Two orthogonal notions of complexity} we will discuss two properties on the forbidden set $U$ that each guarantee that the essential edge number is bounded above \emph{linearly} in terms of the disjoint edge number. In Section \ref{A characterisation in two dimensions} we will prove our main results in the graph case, Theorem \ref{Main theorem in the two-set case} and Proposition \ref{exponential lower bound}. Finally, the more involved proof of our result in the hypergraph case, Theorem \ref{Main theorem for hypergraphs} will unfold in Section \ref{A characterisation for all dimensions}. 

Throughout the remaining sections, all graphs and hypergraphs will be directed, and be allowed to have loops. All hypergraphs will be uniform hypergraphs.

\section{Two guarantees for linear bounds} \label{Two orthogonal notions of complexity}

In this section we discuss several classes of examples of sets $U \subset V \times V$ for which the essential edge number is bounded above \emph{linearly} by the disjoint edge number, uniformly over the choice of graph $G$. We note furthermore that in all examples in this section, this bound will be linear.

\subsection{Rectangle description complexity} \label{Rectangle description complexity} 

The simplest example that we will consider is the case where the complement of $U$ inside $V \times V$, which we will henceforth denote by $U^c$, is a Cartesian product. Although it can appear slightly surprising to consider the complement of $U$ to be this simple set rather than $U$ itself, this is in fact natural to do since the set $U^c$ is the set of edges that are allowed to be in the induced subgraphs arising in the definition of the disjoint edge number of a graph.

\begin{example}\label{Example 1} Let $V$ be a finite set, and let $U$ be a subset of $V \times V$ such that $U^c = X_0 \times Y_0$ for some $X_0,Y_0 \subset V$. Then $e_U(G) = d_U(G)$ for every graph $G=(V,E)$. \end{example}

\begin{proof}

By definition the number $e_U(G)$ is the size of $E \setminus U$, so $e_U(G) = |E \cap (X_0 \times Y_0)|$. Meanwhile, the pairs $X,Y \subset V$ such that $X \times Y$ avoids $U$ are the pairs satisfying $X \times Y \subset X_0 \times Y_0$, which for $X,Y$ both non-empty is equivalent to $X_0 \subset X$ and $Y_0 \subset Y$. The maximum edge number over all such induced subgraphs $e(G[X,Y])$ is $e(G[X_0,Y_0]) = |E \cap (X_0 \times Y_0)|$. The result follows.\end{proof}

We may next consider the situation where rather than being equal to one Cartesian product, the set $U^c$ is equal to a union of a small number of such products.

\begin{example}\label{Example 2}

Let $V$ be a finite set, let $r$ be a positive integer and let $U$ be a subset of $V \times V$ such that \[U^c = (X_1 \times Y_1) \cup \dots \cup (X_r \times Y_r) \]for some $X_1,\dots,X_r,Y_1,\dots,Y_r \subset V$. Then $e_U(G) \le r d_U(G)$ for every graph $G=(V,E)$.

\end{example}

\begin{proof}

By definition the number $e_U(G)$ is the size of the intersection of $E$ with $U^c$. By the union bound this size is at most \[|E \cap (X_1 \times Y_1)| + \dots + |E \cap (X_r \times Y_r)|.\] Meanwhile the sets $X_i \times Y_i$ with $i=1,\dots,r$ are all contained in $U^c$, so \[d_U(G) \ge |E \cap (X_r \times Y_r)|\] for every $i=1,\dots,r$. The result follows.\end{proof}

After that, we turn to the case where $U$ is itself a Cartesian product.

\begin{example}\label{Example 3} Let $V$ be a finite set, and let $U$ be a subset of $V \times V$ such that $U = X_0 \times Y_0$ for some $X_0,Y_0 \subset V$. Then $e_U(G) \le 2d_U(G)$ for every graph $G=(V,E)$. \end{example}

\begin{proof}

We provide two proofs. The first proof consists in partitioning $V \times V$ into four rectangular regions: $U, X_0^c\times Y_0, X_0 \times Y_0^c, X_0^c \times Y_0^c$; letting $e_{2,1},e_{1,2},e_{2,2}$ be the respective sizes of the intersections of $E$ with the last three regions we then note that \[e_U(G) = e_{2,1} + e_{1,2} + e_{2,2}\] whereas \[d_U(G) = \max(e_{2,1} + e_{2,2}, e_{1,2} + e_{2,2}).\] The second proof consists in noticing that from $U = X_0 \times Y_0$ it follows that $U^c$ is the union $(X_0^c\times Y_0) \cup (V \times Y_0^c)$ of two rectangles, reducing to Example \ref{Example 2}.\end{proof}

The second proof will end up being the more easily generalisable of the two, as we next move to the case where $U$ is a union of a bounded number of rectangles.

\begin{example}\label{Example 4}

Let $V$ be a finite set, let $r$ be a positive integer and let $U$ be a subset of $V \times V$ such that \[U = (X_1 \times Y_1) \cup \dots \cup (X_r \times Y_r) \]for some $X_1,\dots,X_r,Y_1,\dots,Y_r \subset V$. Then $e_U(G) \le 2^r d_U(G)$ for every graph $G=(V,E)$.

\end{example}

\begin{proof}

From the decomposition of $U$ we may write \[U^c = \bigcup_{B \subset [r]} \left(\bigcap_{i \in B} X_i^c\right) \times \left(\bigcap_{i \notin B} Y_i^c\right),\] reducing to Example \ref{Example 2} with a union of $2^r$ rectangles. \end{proof}

\subsection{Fibre size complexity}

All examples in the preceding subsection relied in one way or another in $U$ or $U^c$ admitting a decomposition into a bounded number of Cartesian rectangles. Yet we know, from the example of diagonal $U$ that we discussed in the introduction, that there exist sets $U$ for which $e_U$ is linearly controlled by $d_U$, yet cannot be decomposed or have their complement decomposed as a union of a bounded number of rectangles: namely, the diagonal $\Delta$, for which $e_\Delta(G) \le 4 d_\Delta(G)$.

One may ask which property of the diagonal is responsible for this phenomenon. Of course, replacing the diagonal by any set of the type $U = \{(x_1,y_1),\dots,(x_{|V|},y_{|V|})\}$ where $(x_1,\dots,x_{|V|})$ are pairwise distinct and $(y_1,\dots,y_{|V|})$ are pairwise distinct would yield precisely the same result, since that set $U$ becomes the diagonal after a relabelling of $x_1,\dots,x_{|V|}$ and a relabelling of $y_1,\dots,y_{|V|}$. The sets $U$ of the above type are those that intersect each row and each column exactly once, and a natural extension is the class of sets intersecting each row at most $r$ times or every column at most $r$ times for some fixed integer $r$. As we are about to show, within this class of sets $U$ the essential edge number of every graph is linearly bounded above by its disjoint edge number. Furthermore, the proof will be a mild generalisation of the standard probabilistic proof from the case $U = \Delta$.

\begin{example}\label{Example 5}

Let $V$ be a finite set, let $r$ be a positive integer and let $U$ be a subset of $V \times V$ such that the intersection $N_1(y) = \{x \in V: (x,y) \in U\}$ has size at most $r$ for all $y \in V$ or such that the intersection $N_2(x) = \{y \in V: (x,y) \in U\}$ has size at most $r$ for all $x \in V$. Then $e_U(G) \le \frac{(r+1)^{r+1}}{r^r} d_U(G)$ for every graph $G=(V,E)$.

\end{example}

\begin{proof}

By symmetry of the roles played by the rows and columns we may reduce to the case where $N_1(y) = \{x \in V: (x,y) \in U\}$ has size at most $r$ for all $y \in V$, that is, the intersection of $U$ with each column has size at most $r$.

We define subsets $X,Y \subset V$ at random using the following procedure. Let $p$ be a probability which we shall fix later. We first choose $X$ by selecting each element of $V$ independently to belong to $X$ with probability $p$ and to not belong to $X$ with probability $1-p$. After $X$ has been chosen, we deterministically take $Y$ to be the largest possible subset of $V$ such that $X \times Y$ has empty intersection with $U$, that is, we take \[Y = \{y \in V: X \times \{y\} \cap U = \emptyset\}.\]

We next estimate the expected size of the intersection $(X \times Y) \cap E$. If $(x,y) \in E \setminus U$, then the probability that $(x,y) \in X \times Y$ when $X,Y$ are chosen according to the above process is at most $p(1-p)^{|N_1(y)|}$. Since $|N_1(y)| \le r$, we obtain \[\P((x,y) \in X \times Y) \le p(1-p)^r,\] which establishes that the expected size of $(X \times Y) \cap U$ is at most $p(1-p)^r e_U(E)$, and hence in particular that \[d_U(H) \ge p(1-p)^r e_U(E).\] Choosing $p = \frac{1}{r+1}$ provides the desired upper bound. \end{proof}

In subsection \ref{Rectangle description complexity} Examples 2.1 - 2.4 could be paired by complements. We may however not do for the above example: as discussed in the introduction, taking $U$ to be the complement of the diagonal as in Example \ref{Complement of diagonal} leads to a disjoint edge number at most $1$ while the essential edge number may be unbounded.

\section{A characterisation in two dimensions} \label{A characterisation in two dimensions}

\subsection{Upper bound on essential edge number in terms of disjoint edge number}

In this subsection we prove Theorem \ref{Main theorem in the two-set case}.

\begin{proof}[Proof of Theorem \ref{Main theorem in the two-set case}]

Let $V$ be a finite set, let $U$ be a finite subset of $V \times V$, let $r$ be the isolation number of $U$, and let $G=(V,E)$ be a directed graph (possibly with loops). Our task is to prove the upper bound \[e_U(G) \le R(r+1, d_U(G)+1).\] We construct an auxiliary undirected loopless graph $\hat{G}$ with vertex set $E \setminus U$, and where $(x,y)$,$(x’,y’)$ are joined by an edge if and only if $(x’,y) \in U$ or $(x,y’) \in U$.

It follows straightforwardly from this definition that the cliques of $\hat{G}$ are the isolated subsets of $E \setminus U$ with respect to $U$.

Meanwhile, the independent sets $\{(x^1,y^1),\dots,(x^t,y^t)\}$ of $\hat{G}$ correspond to the pairs $(X,Y)$ of subsets of $V$ satisfying $(X \times Y) \cap U = \emptyset$, through \[X = \{x^1,\dots,x^t\}, Y = \{y^1,\dots,y^t\}.\]

We next consider a coloured complete graph $K(\hat{G})$ with vertex set $V(\hat{G}) = E \setminus U$, where an edge is coloured in red if it is an edge of $\hat{G}$ and in blue if it is not an edge of $\hat{G}$.

If the size $e_U(G)$ of $|E \setminus U|$ is at least $R(r+1, d_U(G)+1)$ then $K(\hat{G})$ contains a red clique with size $r+1$ or a blue clique with size $d_U(G)+1$, so $\hat{G}$ contains a clique set with size $r+1$ or an independent set with size $d_U(G)+1$, and therefore $U$ contains an isolated set with respect to $U$ with size $r+1$, or there exists $X \times Y$ with size at least $d_U(G)+1$ such that $(X \times Y) \cap U = \emptyset$. Both possibilities lead to an immediate contradiction.\end{proof}

\subsection{Linear-exponential lower bound in terms of isolation number}

In this subsection we prove Proposition \ref{exponential lower bound}, establishing that the exponential bound in Theorem \ref{Main theorem in the two-set case} is tight up to the constant in front of the isolation number in the exponent.

\begin{proof}[Proof of Proposition \ref{exponential lower bound}]

Let $n=|V|$. We take $U$ at random as follows: the diagonal elements $(x,x) \in V^2$ deterministically do not belong to $U$, and every element $(x,y)$ with $x \neq y$ belongs to $U$ with probability $p$, for some $p \in (0,1)$ that we shall fix later. We take the edge set $E$ of $G$ to be the set of diagonal elements $(x,x) \in V^2$, so we have deterministically $e_U(G)=n$.

We next bound the isolation number $r$ of $U$ from above. For an integer $r_0$, the isolation number of $U$ is at least $r_0+1$ if and only if there exist $(x^1,y^1),\dots,(x^{r_0+1},y^{r_0+1}) \in U^c$ such that for every $1 \le i < j \le r_0+1$ we have $(x^i,y^j) \in U$ or $(x^j,y^i) \in U$. For a given set $r_0+1$-tuple of points as above, the probability that this second condition is satisfied is at most $(2p)^{\binom{r_0+1}{2}}$; as the number of such $r_0+1$-tuples is at most $n^{2(r_0+1)}$, by the union bound we obtain that \[\P(r > r_0) \le n^{2(r_0+1)} (2p)^{\binom{r_0+1}{2}} \le n^{2(r_0+1)} (2p)^{r_0(r_0+1)/2}.\] We then choose $p =\frac{1}{4} n^{-4/r_0}$, so that we obtain \[\P(r > r_0) \le 2^{-r_0(r_0+1)} \le \frac{1}{4}.\]

Afterwards, we bound $d_U(G)$ from above. Due to our choice of edge set of $G$, the number $d_U(G)$ is the size of the largest set $X \subset V$ such that $X \times X$ has empty intersection with $U$. For an integer $d_0$, we obtain \[\P(d_U(G) \ge d_0) \le \binom{n}{d_0} (1-p)^{\frac{d_0(d_0-1)}{2}} \le e^{d_0 \log N - pd_0(d_0-1)}.\] We take $d_0 = \lceil 3p^{-1} \log n \rceil$, which tends to $0$ sufficiently quickly for the event $d \le d_0$ to be useful for us while at the same time sufficiently slowly so that the upper bound on its probability to tend to $0$.

We have established that for $n$ large enough, the above construction provides with positive probability that simultaneously $e_U(G) = n$, the isolation number of $U$ is at most $r$, and $d_U(G) \le 12 n^{4/r} \log n$. This concludes the proof. \end{proof}

\section{A characterisation for all dimensions} \label{A characterisation for all dimensions}

In proving Theorem \ref{Main theorem for hypergraphs}, we can assume without loss of generality that the set $E$ of (hyper)edges of $G=(V,E)$ has empty intersection with $U$ (as replacing $E$ by $E \setminus U$ leaves $e_U(G)$ and $d_U(G)$ unchanged) and will do so throughout the proof. It then also follows immediately that $e_U(G) =|E|$.

\subsection{A tree construction}

If $A$ is a subset of $V^d$, we can choose a subset $M(A) \subset A$ of maximal cardinality with enclosure that avoids $U$, that is, $\mathcal{E}(M(A)) \subset U^c$ but $\mathcal{E}(M(A) \cup \{x\}) \cap U \neq \emptyset$ for every $x \in A \setminus M(A)$. It follows immediately from this definition that $M(A)$ has size $d_U(A)$. (In particular, if $A$ is not contained in $U$, then $M(A)$ is non-empty.)

If $x,y \in V^d$, and $J \subset [d]$ is a set of coordinates, then we write $m_J(x,y)$ for the point of $V^d$ with coordinates in $J$ equal to those of $x$ and coordinates in $J^c$ equal to those of $y$: more formally, if $m= m_J(x,y)$ then $m_j = x_j$ for every $j \in J$ and $m_j = y_j$ for every $j \in J^c$. If $x \in V^d$ and $A \subset V^d$, then we write $A_{J,x}$ for the set of $y \in A$ such that $m_{J}(x,y) \in U$.

Our first lemma informally states that a set $E$ of hyperedges of a $d$-uniform ordered hypergraph that avoids $U$ can always be partitioned into $M(E)$ and a covering of $E \setminus M(E)$ using the sets $E_{J,x}$ with $x \in \mathcal{E}(M(E))$.

\begin{lemma} \label{covering lemma}

Let $d \ge 2$ be an integer, let $V$ be a finite set, let $H=(V,E)$ be a $d$-uniform ordered hypergraph, and let $U$ be a subset of $V^d$ such that the intersection $E \cap U$ is empty. Then we have \[E \setminus M(E) = \bigcup_{x \in \mathcal{E}(M(E))} \bigcup_{J \subset [d]: J \neq \emptyset, [d]} E_{J,x}.\]

\end{lemma}

\begin{proof}

We first note that the sets $E_{J,x}$ are contained in $E \setminus M(E)$. Indeed, if $x \in \mathcal{E}(M(E))$ and $y \in M(E)$, then for every $J \subset [d]$ we have that $m_J(x,y)$ belongs to the enclosure $\mathcal{E}(M(E))$, which has empty intersection with $U$, so $m_J(x,y) \in U^c$.

There remains to check that the sets $E_{J,x}$ cover $E \setminus M(E)$. If $y \in E \setminus M(E)$, then by maximality of $M(E)$ the enclosure $\mathcal{E}(M(E) \cup \{y\})$ has non-empty intersection with $U$. Let $u$ be a point of that intersection. We write \[\mathcal{E}(M(E)) = L_1 \times \dots \times L_d.\]

Let $J$ be the set of coordinates $j \in [d]$ such that $u_j \in L_j$. We cannot have $J = [d]$: if that were the case, then we would have $u \in \mathcal{E}(M(E))$, which is excluded since $u \in U$ but $\mathcal{E}(M(E))$ is contained in $U^c$. We also cannot have $J = \emptyset$: if that were the case, then we would have $u = y$, which is excluded since $u \in U$ but $y$ is contained in $E$, which is a subset of $U^c$. Letting $x$ be an element of $\mathcal{E}(M(E))$ such that $x_j = u_j$ for every $j \in J$ (which exists, since $M(E)$ is not empty), we obtain $y \in E_{J,x}$. \end{proof}

Lemma \ref{covering lemma} then allows us to inductively define a tree, as follows. We take the root to be labelled by $E$, and all other vertices of the tree to be labelled by subsets of $E$. At the root, the immediate descendants are labelled by the sets $E_{J,x}$ with $x \in \mathcal{E}(M(E))$ and $J \subset [d]$ with $J \neq \emptyset, [d]$, unless $E= M(E)$, in which case the tree stops at the root. More generally, if at any point in the tree we have a vertex labelled some $E’ \subset E$ such that $E’\neq M(E’)$, then we define the descendants of $E’$ in the same manner: they are the sets $E_{J,x}$ with x$ \in \mathcal{E}(M(E))$ and $J \subset [d]$ with $J \neq \emptyset, [d]$. We note that if $E’$ is an immediate descendant of $E$, then $|E’|<|E|$, which ensures that there cannot be infinite branches. The leaves are then the vertices labelled $E’$ with $E’=M(E’)$.

From this construction it follows that the set $E$ is the union of the sets $M(E’)$ over all vertices $E’$ of the tree. Since each vertex $E’$ is a subset of $E$, every subset of $E’$ with enclosure avoiding $U$ is also a subset of $E$ with enclosure avoiding $U$, so we have $|M(E’)| \le |M(E)|$; since $|M(E)| = d_U(H)$, it follows that $|M(E’)| \le d_U(H)$. The union bound therefore shows that the size of $E$ is at most $d_U(H)$ times larger than the number of vertices of the tree. Throughout the remainder of the proof, our focus will hence be on obtaining a reasonable upper bound on the size of the tree.

\subsection{Bounds on the tree size}

The next lemma provides a bound on the maximum number of immediate descendants of any vertex of the tree.

\begin{lemma} \label{degree bound}

Every vertex of the tree has at most $2^d d_U(H)^d$ immediate descendants.

\end{lemma}

\begin{proof}

We consider a vertex of the tree labelled $E'$. The number of sets $J \subset [d]$ with $J \neq \emptyset, [d]$ is at most $2^d-2 \le 2^d$, and the size of $\mathcal{E}(M(E’))$ is at most $|M(E’)|^d$. As discussed above, since $E’$ is a subset of $E$, we have $|M(E’)| \le |M(E)|$. In turn we have $|M(E)| = d_U(H)$, so the size of $\mathcal{E}(M(E’))$ is hence at most $d_U(H)^d$. Therefore, the number of immediate descendants of $E$ is at most $2^d d_U(H)^d$.\end{proof}

We next turn to the maximum depth of the tree. Since the labels of vertices on a branch decrease strictly (with respect to inclusion) with the depth, an upper bound is $|E|$, which is not satisfactory for our purposes, since our goal is to obtain an upper bound on $|E|$ using the tree. The construction of the sets $E_{J,x}$ and of the tree nonetheless forces the depth to be bounded in terms of $d_U(H)$ and the isolation number of $E$ (as well as $d$). Together with Lemma \ref{degree bound}, the following lemma will suffice to deduce such an upper bound on the $|E|$ (and hence on $e_U(E)$).

\begin{lemma} \label{depth bound}

If $U$ has isolation number $r$, then the tree has depth at most $2^dr$.

\end{lemma}

\begin{proof}

Let $s$ be the maximum length of any branch of the tree. For every $i=0,1,\dots,s$, let $E_i$ be the label of the depth-$i$ vertex of the branch (with in particular $E_0=E$). For each $i \le s-1$ let $J_i \subset [d]$ such that $E_{i+1} = (E_i)_{J_i,x_i}$ for some $x_i \in \mathcal{E}M((E_i))$. The number of possibilities for $J_i$ is at most $2^d$, so by the pigeonhole principle there exist $t \ge s/2^d$ and $0 \le i_1 < \dots < i_t \le s-1$ such that the sets $J_{i_1},\dots,J_{i_t}$ are all equal to the same set $J$.

For each $h \in [t]$ we choose an element $c_{i_h}$ of $M(E_{i_h})$ and define \[p_{i_h} = m_{J,x_{i_h}}(c_{i_h}).\] What we will show is that the elements $p_{i_h}$ with $h \in [t]$ constitute an isolated set.

Let $h \in [t]$. Because the coordinates of $p_{i_h}$ all come from $x_{i_h}$ and $c_{i_h}$, which both belong to the enclosure $\mathcal{E}(M(E_{i_h}))$, the element $p_{i_h}$ also belongs to this enclosure, and in particular does not belong to $U$.

There remains to establish, using the process through which the points $p_{i_h}$ were constructed, that the enclosure of any two distinct such points has non-empty intersection with $U$. Let $1 \le h’<h'' \le t$. Since $i_{h''} \ge i_{h’}+1$, the set $E_{i_{h''}}$ is a subset of $E_{i_{h’+1}}$. It then follows from the definition \[E_{i_{h’+1}} = \{y \in E_{i_{h’}}: m_{J}(x_{i_{h’}},y) \in U \}\] that $m_{J}(x_{i_{h’}}, c_{i_{h''}})$ belongs to $U$. Let us hence call this point $u_{h’,h''}$. By construction, the coordinates of $u_{h’,h''}$ can be partitioned into coordinates inherited from $p_{i_{h’}}$ and coordinates inherited from $p_{i_{h''}}$, so $u_{h’,h''}$ belongs to the enclosure $\mathcal{E}(\{p_{i_{h’}}, p_{i_{h''}}\})$.

Because the isolation number of $U$ is equal to $r$, we have $t \le r$, so $s \le 2^d r$. \end{proof}

We finally conclude the argument by upper bounding the number of vertices in the tree in the straightforward manner using the previous two lemmas.

\begin{proof}[Proof of Theorem \ref{Main theorem for hypergraphs}]

Each vertex of the tree has at most $2^d d_U(H)^d$ immediate descendants, and the tree has depth at most $2^dr$. The number of vertices of the tree is hence at most \[\sum_{i=0}^{2^dr} (2^d d_U(H)^d)^i \le 2 (2^d d_U(H)^d)^{2^dr}.\] Therefore, we conclude \[e_U(H) = |E| \le 2d_U(H) (2^d d_U(H)^d)^{2^dr} = (2d_U(H))^{d2^dr+1}. \qedhere\] \end{proof}

\end{document}